\documentclass[12pt]{amsart}

\newcounter{theorem}
\makeatletter
\@addtoreset{equation}{section}
\@addtoreset{theorem}{section}
\makeatother

\def\vs{\vskip 1pc}

\renewcommand{\thetheorem}{\thesection.\arabic{theorem}}

\newenvironment{theorem}
{\refstepcounter{theorem} \vs\noindent \bf  \thetheorem.
Theorem. \  \it  }

\newenvironment{lemma}
{\refstepcounter{theorem}\vs\bf\noindent  \thetheorem. Lemma. \ \it}

\newenvironment{prop}
{\refstepcounter{theorem}\vs\bf\noindent  \thetheorem.
Proposition. \ \it}

\newenvironment{definition}
{\refstepcounter{theorem}\vs\bf \noindent
\thetheorem. Definition. \ \rm}

\newenvironment{remark}
{\refstepcounter{theorem}\vs\bf \noindent \thetheorem. Remark. \ \rm}

\newenvironment{corollary}
{\refstepcounter{theorem}\vs\bf\noindent  \thetheorem.
Corollary. \ \it}

\begin{document}

\def\C{{\mathbf C}}
\def\R{{\mathbb R}}
\def\H{{\mathbf H}}
\def\Pt{\tilde{\mathbb P}}
\def\P{{\mathbb P}}
\def\E{{\mathbb E}}
\def\Et{\tilde{\mathbb E}}
\def\Pas{{\mathbf P}{\rm - a.s.}}
\def\bS{{\mathbf S}}
\def\bu{{\mathbf u}}
\def\bU{{\mathbf U}}
\def\L{{\mathcal L}}
\def\A{{\mathcal A}}
\def\B{{\mathcal B}}
\def\cF{{\mathcal F}}
\def\cH{{\mathcal H}}
\def\cD{{\mathcal D}}
\def\cM{{\mathcal M}}
\def\PP{{\mathcal P}}
\def\N{{\mathcal N}}

\def\Om{\Omega}
\def\a{\alpha}
\def\sg{\sigma}
\def\ep{\varepsilon}
\def\dl{\delta}
\def\D{\Delta}
\def\ld{\lambda}
\def\Ld{\Lambda}
\def\gm{\gamma}
\def\Gm{\Gamma}
\def\k{\kappa}
\def\bt{\beta}
\def\f{\varphi}
\def\F{\Phi}


\def\cd{\cdot}
\def\to{\rightarrow}
\def\pd{\partial}
\def\sumr{\sum_{l=1}^r}
\def\Cb{\bar{C}}
 \def\as{\asymp}
\def\FF{I\!\!F}
\def\uth{u^{\theta}}
\def\ut0{u^{\theta^0}}
\def\th{\theta}
\def\tht{\hat{\theta}}

\def\vsm{\vskip -0.5pc}

\def\endp{\ \ \hfill $\Box$ }

\def\ni{\noindent}
\def\nl{\newline}
\def\npb{\nopagebreak}

\def\be{\begin{equation}}
\def\ee{\end{equation}}

\def\ds{\displaystyle}

\def\ba{\begin{array}}
\def\ea{\end{array}}

\def\bd{\begin{description}}
\def\ed{\end{description}}

\def\bee{$$}

\def\bdf{\begin{definition}}
\def\edf{\end{definition}}

\def\br{\begin{remark}}
\def\er{\end{remark}}

\def\bc{\begin{corollary}}
\def\ec{\end{corollary}}

\def\bt{\begin{theorem}}
\def\et{\end{theorem}}

\def\bl{\begin{lemma}}
\def\el{\end{lemma}}

\def\pr{{\bf Proof. }}

\title{Nonlinear Filtering of Diffusion Processes in
Correlated Noise:
Analysis by Separation of Variables}
\author{Sergey V.  Lototsky\\
{$  \ \ \       $ } \\
Department of Mathematics\\ University of Southern California\\
 1042 Downey Way\\ Los Angeles, CA 90089-1113\\
Tel. (213) 740--2389, Fax (213) 740--2424\\
e-mail: {\tt lototsky@math.usc.edu}}

\begin{abstract}
An approximation to the solution of a stochastic parabolic equation
is constructed using the Galerkin approximation followed by the
Wiener Chaos decomposition. The result is applied to
the nonlinear filtering problem  for the
time homogeneous diffusion model with correlated noise. An algorithm is
proposed for computing  recursive approximations of the unnormalized
filtering density and filter, and the errors of the approximations are estimated.
Unlike most
existing algorithms for nonlinear filtering, the real-time part of the
algorithm does not require solving partial differential equations or
evaluating integrals. The algorithm can be used for both continuous and
discrete time observations.
\end{abstract}
\maketitle

\vskip 0.2in

{\bf Short title:} Wiener chaos for nonlinear filtering.

\section{Introduction}

Let $(\Omega, {\mathcal F}, {\mathbb P})$ be a probability space
 with a Wiener process $W=W(t)$.
 Consider a random field $u=u(t,x)$, $0\leq t\leq T$, $x\in \R^d$ so that
$u\in L_2(\Omega; \C([0,T]; L_2(\R^d)))$ and $u(T, \cdot)$  is measurable with
respect to the sigma-algebra $\cF^W_T$ generated by the Wiener process
up to time $T$. If $\{ e_k, k\geq 1\}$ is an orthonormal basis in $L_2(\R^d)$ and
$\{\xi_m, m\geq 1\}$ is an orthonormal basis in $L_2(\Omega,\cF^W_T)$,
then we can write
\begin{equation}
\label{eq:intr1}
u(t,x)\approx \sum_{k=1}^K\left( \sum_{m=1}^M \varphi^K_{m,k}(t)\xi_m
\right)e_k(x).
\end{equation}
where $\varphi_{km}(T)$ are some deterministic coefficients. The objective of the
current work is to study approximation (\ref{eq:intr1}) when the random field
$u$ is a solution
of a stochastic parabolic equation. For such random fields it is possible to
derive explicit representation for the coefficients $\varphi_{km}$ and get an
upper bound on the approximation error in (\ref{eq:intr1}). The results are then
used to derive an approximate algorithm for solving the nonlinear filtering
problem of diffusion processes with correlated noise.

 The problem of nonlinear filtering can be briefly described  as follows.
 Assume that $(X=X(t), Y=Y(t)), \ t  \geq 0,$ are two diffusion processes
 with values in $\R^d$ and $\R^r$ respectively,  so that
 $X$ is the unobservable component and the observable component
 $Y$ is given by
\bee
Y(t)=\int_0^t h(X(s))ds + W(t).
\bee
The problem is called noise-uncorrelated, if the Wiener process  $W=W(t)$,
 representing the observation noise, is independent of $X$. The problem is called
 noise-correlated,  if there is correlation between $W$ and $X$.
  If $f=f(x)$ is a measurable function satisfying
$\E |f(X(t))|^2 < \infty, \ t\geq 0,$ then the problem of nonlinear filtering is
to find the best mean square estimate
$\hat{f}_t$ of $f(X(t))$ given the trajectory $Y(s),\  s \leq t.$
It is known \cite{Kal,LSh.s,Roz}
that, under certain regularity assumptions, we have
\be
\label{eq:intr3}
\hat{f}_t=\frac{\int_{\R^d} f(x)p(t,x)dx}{\int_{\R^d} p(t,x)dx},
\ee
where $p=p(t,x)$ is a random field called the unnormalized filtering
density (UFD). The problem of estimating $f(X(t))$ is thus reduced to the
problem of computing the UFD $p$. It is also known \cite{Roz} that $p=p(t,x)$ is
the solution of the Zakai filtering equation, a stochastic parabolic  equation,
 driven by the observation process.
The exact solution of this equation can be found only in some
special cases, and the development of numerical schemes for solving
the Zakai equation has become an area of active research.

Many  of the existing numerical schemes for the Zakai equation
use various generalizations of the corresponding algorithms for
 the deterministic partial differential equations.  Examples of the corresponding
algorithms
can be found in Bennaton \cite{Bennaton}, Florchinger and LeGland
\cite{FlLg}, Ito \cite{KIto}, etc. Because of the
 large amount of calculations, these algorithms
cannot be implemented in real time when the dimension of the
state process is more than three.

In some applications, like target tracking, the filter estimate must
be computed in real time. Such applications require filtering
algorithms with fast on line computations.
When the parameters of the model are known in advance,
 the real time computations can be simplified by  separating
  the deterministic and stochastic
components of the Zakai equation and performing
the computations related to the deterministic component
in advance. The separation is based on the Wiener chaos decomposition
of solutions of stochastic parabolic equations.
Starting with the works of Kunita \cite{Kunita81}, Ocone \cite{Ocn}, and
Lo and Ng \cite{LoNg}, this approach was further developed by
Budhiraja and
Kallianpur \cite{BdKl3,BdKl1,BdKl2} and
Mikulevicius and Rozovskii \cite{MR.sep,MR.fhe,mr2,mir2}.
An  algorithm to solve the Zakai
equation using this approach for the noise uncorrelated problem
 was suggested in Lototsky et al. \cite{LMR}.
The algorithm in \cite{LMR} was based on the following approximation.
First, the unnormalized filtering density was approximated using the
Wiener Chaos decomposition:
\be
\label{eq:intr2}
p(t,x)\approx \sum_{m=1}^M \varphi_m(t,x)\xi_m.
\ee
After that, the coefficients $\varphi_m$ were expanded in the basis
$\{e_k\}$ in $L_2(\R^d)$, resulting in the approximation
\be
\label{eq:intr22}
p(t,x)\approx \sum_{m=1}^M\left( \sum_{k=1}^K \varphi_{m,k}(t)e_k(x)
\right)\xi_m.
\ee
In other words, first, the stochastic variable
was separated, and then, the spacial variable.

Alternatively, one can start with
 the Galerkin approximation of $p$:
\be
\label{eq:intr5}
p(t,x)\approx \sum_{k=1}^K p^K_k(t)e_k(x).
\ee
The coefficients $p^K_k(t)$ satisfy a system of stochastic ordinary
differential equations driven by the observation process.
The solution of this   system
is then expanded   using the  Wiener Chaos
decomposition, resulting in the approximation of the type (\ref{eq:intr1}):
\be
\label{eq:intr55}
p(t,x)\approx \sum_{k=1}^K\left( \sum_{m=1}^M \varphi^K_{m,k}(t)\xi_m
\right)e_k(x).
\ee
 In other words, first, the spacial variable is separated,
then, the stochastic variable. The optimal filter is then approximated by
\be
\label{eq:intr555}
\hat{f}_t\approx \frac{\sum_{k=1}^K\left( \sum_{m=1}^M
\varphi^K_{m,k}(t)\xi_m
\right)f_k}{\sum_{k=1}^K\left( \sum_{m=1}^M \varphi^K_{m,k}(t)\xi_m
\right)\alpha_k},
\ee
where $f_k=\int_{\R^d}f(x)e_k(x)dx, \alpha_k=\int_{\R^d}e_k(x)dx$.
First suggested in \cite{CPF} as a computational alternative to
(\ref{eq:intr22}), approximation (\ref{eq:intr55}) was further
analysed in \cite{LFung}.

The order in
which the variables are separated does make a difference. The algorithms
based on approximations (\ref{eq:intr22}) and (\ref{eq:intr55}) have
different approximation errors and, unlike (\ref{eq:intr22}),
analysis of (\ref{eq:intr55}) is possible for noise correlated problem.

Recall that the Zakai filtering equation for the unnormalized filtering
density $p=p(t,x)$ is
\be
\label{eq:intr_ze}
dp={\mathcal L}^* p\;dt + {\mathcal M}^*p \;dY(t).
\ee
The elliptic differential operator ${\mathcal L}$ is the generator of
the unobserved process $X$, while the operator ${\mathcal M}$
is bounded in the noise uncorrelated problem and in unbounded in the
noise correlated problem. The presence of the unbounded
operator in the stochastic part of  equation (\ref{eq:intr_ze})
for  the noise correlated problem
 makes the analysis and implementation of the numerical methods
 for the Zakai equation much more difficult  (see, for
 example, Florchinger and LeGland
\cite{FlLg}).

The objective of the current work is to analyze the   algorithm
for solving the Zakai equation using  approximation (\ref{eq:intr55}).
 First, (\ref{eq:intr55}) is studied for an abstract stochastic evolution
 system. In Section 2, the Galerkin approximation is investigated,
 and in Section 3, the Wiener chaos decomposition for a system
 of stochastic ordinary differential equations. In each situation,
 the rate of convergence is established in terms of the numbers
 $K$ and $M$ of the basis function used. The filtering problem
 is introduced in Section 4, the filtering algorithm is presented in
 Section 5, and the convergence of the algorithm is studied in
 Section 6. The real time part of the proposed algorithm does not
 require solving  differential equations or using quadrature
 methods to evaluate integrals in (\ref{eq:intr3}). The algorithm can
 also be used if the observations are available in discrete time.

 Unlike the previous works on the subject, this paper presents a
 unified treatment of both noise-correlated and noise-uncorrelated
 problems with possibly degenerate diffusion in the un-observed
 component. Another difference from the previous works on the
 subject is that the error bound is derived not only for the
 filtering density but also for the optimal filter $\hat{f}_t$
 with a large class of functions $f$.

\section{Galerkin approximation of stochastic evolution equations}
\setcounter{equation}{0}

Consider the stochastic evolution system
\be
\label{es33}
u(t) = u_0+\int_{T_0}^t  \A u(s)ds +
\int_{T_0}^t\sumr \B_l u(s) dW_l(s),\ T_0\geq 0,
\ee
 where
 $ W=W(t)$ is an $r$ - dimensional standard Wiener process on a
 complete
probability space $(\Om, \cF, \P)$, $u_0$ is independent of $W$,  and
 $\A$ and $\B_l,\ l=1,\ldots, r,$ are  linear operators
 acting in the scale of infinite dimensional
 Hilbert  spaces $\{\H^a, a \in \R\}$.
   To simplify the notation, both the inner product in $\H^0$ and
the duality between $\H^1$ and $\H^{-1}$ will be denoted by
$(\cdot, \cdot)_0$; $\|\cdot\|_a$ is the norm in the space $\H^a$.
 It will be assumed that  equation
(\ref{es33}) is either coercive or dissipative \cite[Chapter 3]{Roz}.
In particular, there exists  a constant $C^*>0$
so that, for every $v\in \H^1$,
\be
\label{eq:diss00}
\|\A v\|_{-1}\leq C^*\|v\|_{1},\ \|\B_lv\|_{0} \leq C^*\|v\|_{1},\
{\rm and}\  2(\A v,v)_0+\sumr \|\B_l v\|^2_0 \leq C^*\|v\|_0^2.
\ee
If $u_0 \in \H^1$, then there is  a unique solution
$u=u(t)$ in the space
 $L_2(\Om\times[T_0,T];\H^1)\cap L_2(\Om;\C([T_0,T];\H^0))$
 (see Theorems 3.1.4 and 3.2.2 in \cite{Roz}).

Suppose there exists an orthonormal basis  $\{e_k,k\geq \}$  in $\H^0$
so  that
$e_k\in \H^1$ for all $k$.  Consider the following system
of stochastic ordinary differential equations:
\be
\label{ga}
\ba{lll}
\ds \ds du^K_k(t)&=& \ds\sum_{n=1}^K ( \A e_n,e_k )_{0}
 u^K_n(t)dt\\
\ds   &+&\ds \sumr \sum_{n=1}^K  (\B_le_n, e_k)_0
u^K_n(t) dW_l(t),\ T_0< t \leq T, \\
\ds u^K_k(T_0) &=& (u_0,e_k)_0,\ k=0, \ldots, K.
\ea
\ee
 The function
\bee
u^K(t)=\sum_{k=1}^K u_k^K(t)e_k
\bee
is called the {\bf Galerkin approximation}  of $u(t)$.
It is proved in the following theorem that, under some natural assumptions,
\bee
\lim_{K \to \infty} \sup_{T_0 \leq t \leq T}
\E \| u(t;T_0;u_0)-u^K(t)\|_0^2 = 0,
\bee
and the rate of convergence is determined.

\bt
\label{H2.th}
Let the following conditions be fulfilled:
\begin{enumerate}
\item[{\bf 1.}]  The basis $\{ e_k \}$ consists of the
eigenfunctions of
a linear operator  $\Ld$
with the corresponding eigenvalues $\ld_k$. The operator $\Ld$ is a symmetric
operator in $\cH^0$ and
 there exist
numbers $0<c_1<c_2$ and $\theta>0$ so that, for all $k$,
\begin{equation}
\label{eq:eigass}
 c_1\leq \ld_kk^{-\theta}\leq c_2;
 \end{equation}
\item[{\bf 2.}] $e_k \in \H^1$ and $\|e_k\|_1 \leq C_ek^q,\ q \geq 0$;
\item[{\bf 3.}] $\sup_{T_0 \leq t \leq T}
 \E \|\Ld^{\nu} u(t) \|_0^2 < \infty$
for some positive integer $\nu$ so that $\theta_1 := \nu \theta-2q>1$.
\end{enumerate}
Then
\be
\label{grc}
\sup_{ T_0 \leq t \leq T}\E \| u(t)-u^K(t)\|_0^2 \leq
\sup_{T_0 \leq t \leq T} \E \|\Ld^{\nu} u(t) \|_0^2 \;
\frac{Ce^{C(T-T_0)}}{K^{2(\theta_1-1)}},
\ee
where  $C$ is a constant depending only on the
constant $C^*$ in (\ref{eq:diss00}) and  the numbers
$c_1, c_2, C_e, \nu,\ \theta,\  q$.
\et

\ni \pr
If $\psi_k(t) := (u(t), e_k)_0$, then
\be
\label{pg1}
\E \|u(t)-u^K(t)\|_0^2 = \sum_{k=0}^K
 \E |\psi_k(t)-u_k^K(t)|^2+
\sum_{k > K} \E |\psi_k(t)|^2.
\ee
By assumptions 1 and  3 of the theorem,
\be
\label{H2.p1}
|\psi_k(t)| \leq \frac{\|\Ld^{\nu}u(t)\|_0}{\ld_k^{\nu}}
\ee
so that
\be
\label{pg2}
\ds \sup_{T_0\leq t \leq T} \sum_{k > K} \E |\psi_k(t)|^2 \ds \leq
\sup_{T_0\leq t\leq T}\E\|\Ld^{\nu}u(t)\|_0^2
\frac{C}{K^{2\nu\theta-1}}\leq
\ds \sup_{T_0\leq t\leq T}\E\|\Ld^{\nu}u(t)\|_0^2
\frac{C}{K^{2(\theta_1-1)}}.
\ee
For $1 \leq k \leq K$ define $\dl_{k}(t):=\psi_k(t)-u_k^K(t),$
so that $\sum_{k=1}^K \E |\psi_k(t)-u^K_k(t)|^2 =
\sum_{k=1}^K \E|\dl_k|^2,$
and also define
\bee
\dl_{1,n}(t):=\sum_{k>K}(\A e_{k},e_{n})_0
 \psi_k(t),\quad
 \dl_{2,n}^l(t):=\sum_{k > K}(\B_le_{k}, e_{n})_0 \psi_k(t).
\bee
Both $\dl_{1,n}(t)$ and $\dl_{2,n}^l(t)$ are well defined due to
 (\ref{H2.p1})
and  assumptions 2 and 3 of the theorem.
Then
\be
\label{pg3}
\ba{lll}
\ds d\dl_{n}(t) &=& \ds \sum_{k=1}^K(\A e_{k},e_{n})_0
 \dl_{k}(t) dt+
\sumr \sum_{k=1}^K(\B_le_{k}, e_{n})_0\dl_{k}(t) dW_l(t) \\
{}&+& \ds\dl_{1,n}(t) dt + \sumr \dl_{2,n}^l dW_l(t), \ T_0 < t \leq T; \\
\ds  \dl_{n}(T_0) &=& 0,\ 1 \leq n \leq K,
\ea
\ee
and by the Ito formula,
\be
\label{pg4}
\ba{l}
\ds  \sum_{n=1}^K \E|\dl_{n}(t)|^2 = 2 \int_{T_0}^t \sum_{n,k =1}^K
 (\A e_{k},e_{n})_0  \E \dl_{n}(s) \dl_{k}(s) ds \\
+ \ds \sumr \sum_{n=1}^K  \int_{T_0}^t
\E\left( \sum_{k=1}^K\ (\B_le_{k}, e_{n})_0\dl_{k}(s)
\right)^2 ds +
 2 \sum_{n=1}^K  \int_{T_0}^t \E \dl_{1,n}(s)\dl_{n}(s) ds\\
+\ds 2 \sumr \sum_{n,k=1}^K  \int_{T_0}^t (\B_le_{k}, e_n)_0
\E \dl_{2,n}^l(s)
 \dl_{k}(s) ds + \sumr \sum_{n=1}^K  \int_{T_0}^t
  \E (\dl_{2,n}^l(s))^2 ds.
\ea
\ee

 It follows from the third inequality in (\ref{eq:diss00}) that
\be
\label{pg5}
\ba{l}
\ds 2 \int_{T_0}^t \sum_{n,k =1}^K
  (\A e_{k},e_{n} )_0 \E \dl_{n}(s) \dl_{k}(s) ds\\
+\ds \sumr \sum_{n=1}^K \int_{T_0}^t
\E \left( \sum_{k=1}^K (\B_le_{k}, e_{n})_0 \dl_{k}(s)
\right)^2 ds \leq   C \sum_{k=1}^K \int_{T_0}^t \E(\dl_{k}(s))^2 ds.
\ea
\ee
 The first two inequalities in (\ref{eq:diss00}) and   assumption  2 imply
\bee
|( \A e_{k},e_{n})_0|
\leq C\|e_k\|_1\|e_n\|_1 \leq Ck^qn^q,\quad
|(\B_le_k,e_n)_0| \leq C\|e_k\|_1\|e_n\|_0 \leq Ck^q,
\bee
so that by (\ref{H2.p1}),
\bee
|\dl_{1,n}(t)|
\leq Cn^q\frac{\|\Ld^{\nu}u(t)\|_0}{K^{\nu\theta-q-1}},\quad
|\dl_{2,n}^l| \leq
 C\frac{\|\Ld^{\nu}u(t)\|_0}{K^{\nu\theta-q-1}},
\bee
and
\be
\label{pg6}
\ba{l}
\ds  \sum_{n=1}^K \int_{T_0}^T \E (\dl_{1,n}(s))^2 ds +
 \sumr \sum_{n=1}^K \int_{T_0}^T \E(\dl_{2,n}^l(s))^2 ds\\
\leq
\ds  (T-T_0)\sup_{T_0\leq t\leq T}
\E\|\Ld^{\nu}u(t)\|_0^2\;\frac{(r+1)C}{K^{2(\theta_1-1)}}.
\ea
\ee
After that   (\ref{pg4})--(\ref{pg6}) and the obvious inequality
$2|ab| \leq a^2+b^2$
imply
\bee
 \sum_{n=1}^K \E|\dl_{n}(t)|^2 \leq C \sum_{n=1}^K
\int_{T_0}^t \E |\dl_{n}(s)|^2 ds+
(T-T_0)\sup_{T_0\leq t\leq T}
\E\|\Ld^{\nu}u(t)\|_0^2\;\frac{(r+1)C}{K^{2(\theta_1-1)}},
\bee
so that by the Gronwall inequality
\bee
\sup_{T_0\leq t \leq T}
\sum_{n=1}^K \E |\dl_{n}(t)|^2 \leq (T-T_0)\sup_{T_0\leq t\leq T}
\E\|\Ld^{\nu}u(t)\|_0^2\;e^{C(T-T_0)}\;\frac{(r+1)C}{K^{2(\theta_1-1)}}.
\bee
Together with (\ref{pg1}) and (\ref{pg2}), the last inequality implies
(\ref{grc}). Theorem \ref{H2.th} is proved.
\endp

\section{Wiener Chaos Expansion}
\setcounter{equation}{0}

 On a complete probability space $(\Om, \cF, \P)$
  consider a system of  stochastic ordinary differential  equations:
\be
\label{es2}
U(t) = U_0+\int_{T_0}^t  A U(s)ds +
\int_{T_0}^t\sumr B_l U(s) dW_l(s)\ T_0\geq 0,
\ee
where $U(t),U_0\in \R^K,$ $A, B_l \in \R^{K\times K},$ the matrices $A,B_l$ are
deterministic,  and $U_0$ is independent of
the $r$-dimensional Wiener process $W$.
The solution of (\ref{es2}) is denoted by
$U(t;T_0;U_0).$

In what follows,
the Wiener chaos decomposition of $U(t;T_0;U_0)$ will be derived
and the
properties of the decomposition studied.

As the first step, recall the construction of an orthonormal basis
in the space
 $L_2(\Omega, \cF^W_{T_0,t}, \P)$ of square integrable random
 variables that are measurable with respect to the $\sigma$-algebra,
  generated by
  the Wiener process up to time $t$.
 Let $\a$ be an {\bf $r$-dimensional multi-index}, that is,  a  collection
 $\a=(\a_k^l)_{1 \leq l \leq r,\ k \geq 1}$ of nonnegative integers
 such that only finitely many of $\a_k^l$ are different from zero.
The set of all such multi-indices will be denoted by $J$.
For $\a\in J$
define $\a ! := \prod_{k,l}(\a_k^l !)$.

For a fixed  $t^*>T_0$ choose a complete orthonormal system
 $\{ m_k \}=\{ m_k(s)\}_{k \geq 1}$ in $L_2([T_0,t^*])$
and define
\bee
\xi_{k,l}=\int_{T_0}^{t^*} m_k(s)dW_l(s)
\bee
 so that  $\xi_{k,l}$ are
independent Gaussian random variables with zero mean and unit
 variance.

If
\be
\label{her1}
H(x) := (-1)^n e^{x^2/2}\frac{d^n}{dx^n}e^{-x^2/2}
\ee
 is the $n$-th
Hermite polynomial, then the collection
\bee
\left\{ \xi_{\a}(W_{T_0,t^*}) := \prod_{k,l}
\left( \frac{H_{\a_k^l}(\xi_{k,l})}{\sqrt{\a_k^l !}}
 \right),  \quad \a \in J \right\}
\bee
\ni is an orthonormal system in
$L_2(\Omega, \cF^W_{T_0,t^*}, \P).$ A theorem of Cameron
and Martin \cite{CM}
shows that $\{ \xi_{\a}(W_{T_0,t^*}) \}_{\a \in J}$ is actually
a basis in that
space.

\ni \bt
\label{wce.th}
If $\eta \in L_2(\Om,\cF^W_{T_0,t^*},\P),$ then
\be
\label{wce.f}
\eta=\sum_{\a \in J} \E [ \eta \xi_{\a}(W_{T_0,t^*}) ]
\xi_{\a}(W_{T_0,t^*})
\ee
and
\bee
\E|\eta|^2=\sum_{\a \in J} |\E  \eta \xi_{\a}(W_{T_0,t^*}) |^2.
\bee
\et

\ni \pr This theorem is proved in \cite{CM} and \cite{Hida79}.
\endp

\ni \begin{theorem}
\label{Th.sp1}
If   $t^*>T_0$ is fixed,
then, for every $s \in [T_0,t^*]$, the solution $U(s;T_0;U_0)$
 can be written as
\be
\label{u}
U(s;T_0;U_0)=
\sum_{\a\in J}\frac 1{\sqrt{\a !}}\f_{\a}(s;T_0;U_0)\xi_{\a}
(W_{T_0,t^*}),
\ee
and  the following Parseval's equality holds:
\be
\label{pu}
\E | U(s;T_0;U_0) |^2 = \sum_{\a \in J} \frac 1{\a !}
E|\f_{\a}(s;T_0;U_0) |^2.
\ee

The coefficients of the expansion are $\R^K$-vector functions   and
 satisfy the recursive system of deterministic equations
\be
\label{S}
\begin{array}{lll}
\ds \frac{\pd \f_{\a}(s;T_0;U_0)}{\pd s} &=
 \ds A \f_{\a}(s;T_0;U_0)+\sum_{k,l}\a_k^l
 m_k(s) B_l \f_{\a(k,l)}(s;T_0;U_0), \quad T_0<s \leq t^*; \\
\ds \f_{\a}(T_0;T_0;U_0) &= U_01_{ \{ | \a |=0 \} },
\end{array}
\ee
where $  \a = (\a_k^l)_{1 \leq l \leq r,\ k \geq 1}  \in J$ and
 $\a(i,j)$ stands for the  multi-index
 $\tilde{\a} = (\tilde{\a}_k^l)_{1 \leq l \leq r,\ k \geq 1}$ with
\be
\label{al}
\tilde{\a}_k^l=
\begin{cases}
\a_k^l & {\rm \ if\ } k \not= i {\rm \ or \ }l \not= j
{\rm \ or\  both\ } \\
\max(0,\a_i^j-1) &{\rm \ if\ } k=i {\rm \ and\ } l=j.
\end{cases}
\ee

\end{theorem}

\ni {\bf Proof.} Assume first that $U_0=g$ is deterministic;
 the Markov property of the solution of (\ref{es2})
 implies that, once the
 derivation is complete, we can replace $g$ with $U_0$.

If $g$ is deterministic, then
  $ U(s;T_0;g) \in L_2(\Om,\cF^W_{T_0,t^*},\P)$
for $s\leq t^*$, and Theorem \ref{wce.th} implies
(\ref{u}) and (\ref{pu}).

To prove that the coefficients satisfy (\ref{S}), define
\bee
P_t(z)=\exp \Big\{ \int_{T_0}^t
\sum_{l=1}^rm_z^l(s)dW_l(s) -\frac{1}{2}
\int_{T_0}^t \sum_{l=1}^r |m_z^l(s)|^2 ds \Big\},\
 T_0 \leq t \leq t^*,
\bee
where
$ m_z^l=\sum_{k \geq 1}m_k(s)z_k^l$ and
 $\{ z_k^l\},\ l=1, \ldots, r, \ k=1,2,\ldots,$ is
 a sequence of real numbers such that
 $\sum_{k,l} |z_k^l|^2 < \infty$.
Then direct computations show that
\bee
\xi_{\a}(W_{T_0,t^*})=\frac{1}{\sqrt{\a !} }
\frac{\pd^{\a}}{\pd z^{\a} } P_{t^*}(z)\Big|_{z=0}\;,
\bee
where
\bee
\frac{\pd^{\a}}{\pd z^{\a} } = \prod_{k,l} \frac{\pd^{\a^l_k}}
{(\pd z_k^l)^{\a^l_k}}\;,
\bee
and also, that
\bee
\E [\eta \xi_{\a}(W_{T_0,t^*})]=\frac{\pd^{\a}}{\pd z^{\a} }
\E[\eta P_{t^*}(z)]\Big|_{z=0}
\bee
for every $\eta \in L_2(\Om,\cF^W_{T_0,t^*},\P).$
Consequently,
\bee
\ba{c}
\ds \f_{\a}(s;T_0;g)=\ds  \frac{\pd^{\a}}{\pd z^{\a} }
\E [U(s;T_0;g) P_{t^*}(z)]\Big|_{z=0}
\\
=\ds \frac{\pd^{\a}}{\pd z^{\a} }
\E [U(s;T_0;g) P_s(z)]\Big|_{z=0}\;,
\ea
\bee
 where the second equality follows from the martingale property of
$P_s(z)$ on
\nl $(\Omega, \{ \cF_{T_0,t}^W \}_{T_0\leq t \leq t^*}, \P).$
It follows from the definition of $P_s(z)$ that
\bee
dP_s(z) =\sumr m_z^l(s) P_s(z) dW_l(s), \
 T_0 \leq  s \leq t; \ P_{T_0}(z)=1.
\bee
Then (\ref{es2}) and the Ito formula imply that
\bee
\ba{c}
\ds U(s;T_0;g)P_s(z)=\ds g \qquad \qquad \qquad \qquad \qquad \\
+\ds \int_{T_0}^s \Big(
A U(\tau;T_0;g) +
\sum_{l=1}^r B_lU(\tau;T_0;g)\Big) m_z^l(\tau) P_{\tau}(z)d\tau  \\
+\ds \int\limits_{T_0}^s\sumr\Big(
B_lU(\tau;T_0;g)+
U(\tau;T_0;g) m_z^l(s)\Big) P_s(z) dW_l(\tau).
\ea
\bee
Taking the expectation on both sides of the last equality and
setting
\newline $ \f(s,z;T_0;g):=\E U(s;T_0;g)P_s(z)$ results in
\bee
\ds  \f(s,z;T_0;g)=\ds g +\int_{T_0}^s \Big(
 A\f(\tau,z;T_0;g)+
{}\ds \sum_{l=1}^r m_z^l(\tau)
 B_l \f(\tau,z;T_0;g) \Big)d\tau.
\bee
Applying the operator $\ds \frac{1}{\sqrt{\a !}}
\frac{\pd^{\a}}{\pd z^{\a} }$
 and setting $z\!=\!0$ yields that the functions $\f_{\a}(s;T_0;g)$
satisfy  (\ref{S}). Theorem \ref{Th.sp1} is proved.
\endp

For a multi-index $\a\in J$ define
\begin{itemize}
 \item $| \a | := \sum_{l,k} \a_k^l $ (length of $\a$);
\item  $ d(\a) := \max \{ k \geq 1: \, \a_k^l > 0 \
\mbox{for some} \ 1 \,
\leq l \leq r\}$   (order of $\a$).
\end{itemize}
To study the rate of convergence  of the series in  (\ref{u}),
it is necessary
to note that the summation $\sum_{\a \in J}$ is double infinite:
\be
\label{SUM}
\sum_{\a \in J}=\sum_{k=0}^{\infty}\sum_{|\a|=k}
\ee
 and there are infinitely many
multi-indices $\a$ with $|\a|=k>0$.

Define $J_N^n=\{ \a \in J:|\a| \leq N,\ d(\a) \leq n\}$ and then
\be
\label{utrn.1}
U_N^n(s;T_0;U_0)=\sum_{\a\in J_N^n}\frac 1{\sqrt{\a
!}}\f_{\a}(s;T_0;U_0)\xi_{\a}(W_{T_0,t^*}).
\ee
Now  the summation in (\ref{utrn.1}) is over a finite set: if $d(\a) \leq n$,
then there are at most $(nr)^k$ multi-indices $\a$ with $|\a|=k$.

\bt
\label{utrn.cor}
Let the constants $C_0, C_1, C_2$ be such that
$| A v |^2 \leq C_0|v|^2$,
$|e^{At} v|^2\leq e^{C_1t}|v|^2$,
 $|B_l v|^2 \leq C_2| v|^2$ for every vector $v\in \R^K$.
If the  basis  $\{m_k\}$ is the Fourier cosine basis
\be
\label{basis}
\!\!\!\!\!\!\!\!\!\!
 m_1(s)\!=\!\frac 1{\sqrt{t^*-T_0}};\  m_k(s)\!=
 \!\sqrt{\frac{2}{t^*-T_0} }
 \cos \left( \frac{\pi (k-1) (s-T_0)}{t^*-T_0} \right), \, k>1; \
T_0\leq s \leq t^*,
\ee
 then
\be
\label{utrNn}
\ba{c}
\E|U(t^*;T_0;U_0)-U^n_N(t^*;T_0;U_0)|^2 \leq
\ds 2e^{\Cb(t^*-T_0)}
\Bigg( \frac{[C_2r(t^*-T_0)]^{N+1}}{(N+1)!}\\
 \ds +2C_2r\frac{(t^*-T_0)^2}{n}[
 \epsilon(B)+
C_0(1+C_2r(t^*-T_0))(t^*-T_0)]\Bigg)\E|U_0|^2,
\ea
\ee
where $\Cb=C_1+C_2r$ and $0\leq \epsilon(B)\leq 4$;
$\epsilon(B)=0$ if the matrices $B_l$ commute
(in particular, if $r=1$).
\et
This Theorem is proved below in Section 7.

If $t^*-T_0=\Delta$, then (\ref{utrNn}) becomes
\begin{equation}
\label{eq:trNn}
\E|U(t^*;T_0;U_0)-U^n_N(t^*;T_0;U_0)|^2 \leq
e^{C\Delta}\left( \frac{(C\Delta)^{N+1}}{(N+1)!}+
\frac{\Delta^2}{n}(\epsilon(B)+C\Delta)\right)\E|U_0|^2,
\end{equation}
and the constant $C$ depends only on the matrices $A$ and $B_l$ in
(\ref{es2}).

\section{Diffusion Filtering Model}
\setcounter{equation}{0}

Let $(\Om,\cF,\P)$ be a complete probability space with  independent
standard Wiener processes $W=W(t)$ and $V=V(t)$ of dimensions $d_1$
and
$r$ respectively. Let $X_0$ be a random variable independent of $W$
and $V$.
In the {\em diffusion filtering model}, the unobserved
$d$ - dimensional state (or signal) process $X=X(t)$ and the
$r$-dimensional
observation process $Y=Y(t)$ are defined by the  stochastic ordinary
differential equations
\be
\label{fmod.sss2}
\ba{l}
\ds dX(t)=b(X(t))dt+\sg(X(t))dW(t)+\rho(X(t))dV(t),\\
\ds dY(t)=h(X(t))dt+dV(t),\ 0 < t \leq T;\\
\ds X(0)=X_0,\quad Y(0)=0,
\ea
\ee
where $b(x)\in\R^d$, $\sigma(x)\in\R^{d\times d_1}$,
 $\rho(x)\in\R^{d\times r}$,
$h(x)\in \R^r$.

{\bf Assumption  R1}. The functions $\sigma$ and $\rho$  are
$\C^3_b(\R^d)$, that is, bounded and
three times continuously differentiable on $\R^d$ so that all the derivatives
are also bounded; the functions $b$ and $h$ are $\C^2_b(\R)$,
and the
random variable $X_0$ has a density $p_0$.

Under Assumption R1 system (\ref{fmod.sss2}) has a unique strong
solution \cite[Theorems 5.2.5 and 5.2.9]{KarShr}.

If $f=f(x)$ is a scalar measurable function on $\R^d$ so that
 $ \sup_{0 \leq t\leq T}\E|f(X(t))|^2 < \infty$,
then the {\em filtering problem} for (\ref{fmod.sss2}) is to find the
best mean square estimate $\hat{f}_t$ of $f(X(t)),\ t \leq T,$ given
the observations
$Y(s),\ 0< s \leq t$.
 Denote by $\cF^Y_t$ the $\sg$-algebra generated  by
$Y(s),\ 0 \leq s \leq t$. Then the properties of the conditional
expectation imply that the solution of the filtering problem is
\bee
\hat{f}_t=\E\left(f(X(t))|\cF^Y_t \right).
\bee
To derive an alternative representation of $\hat{f}_t$, some
 additional
constructions will be necessary.

Define a new probability measure $\Pt$ on $(\Om,\cF)$ as follows:
for $A \in \cF$,
\bee
\Pt(A)=\int_AZ_T^{-1}d\P,
\bee
where
\bee
Z_t = \exp \left\{ \int_0^t h^*(X(s))dY(s)-\frac{1}{2} \int_0^t
|h(X(s))|^2ds \right\}
\bee
(here and below, if $\zeta \in \R^k$, then $\zeta$ is a {\it
column} vector, $\zeta^*=(\zeta_1, \ldots,\zeta_k),$ and
$|\zeta|^2=\zeta^*\zeta$). If  the function $h$ is bounded, then
the measures $\P$ and $\Pt$ are equivalent.
The expectation with respect to the measure $\Pt$ will be
denoted by $\Et$.

The following properties of the measure $\Pt$ are well known
\cite{Kal,Roz}:
\begin{enumerate}
\item[{\bf P1.}] Under the measure $\Pt$, the distributions of the
Wiener process $W$ and the random variable $X_0$ are
 unchanged, the observation process $Y$ is
a standard Wiener process, and the state process $X$ satisfies
\bee
\ba{l}
\ds
\!\!\!\!\!\!\!\!\!\!
 dX(t)=b(X(t))dt+\sigma(X(t))dW(t)+\rho(X(t))
 \left( dY(t)-h(X(t))dt\right),
\ 0 < t \leq T;\\
\ds \!\!\!\!\!\!\!\!\! X(0)=X_0;
\ea
\bee
\item[{\bf P2.}] Under the measure $\Pt,$  the Wiener processes
$W$ and $Y$
and the random variable $X_0$ are independent of one another;

\item[{\bf P3.}] The optimal filter $\hat{f}_t$ satisfies
\be
\label{bf}
\hat{f}_t=\frac{\Et
\left[ f(X(t))Z_t|\cF_t^Y \right] }{\Et[Z_t|\cF^Y_t]}.
\ee
\end{enumerate}

Because of property {\bf P2} of the measure $\Pt$ the filtering problem will be
studied on the probability space $(\Omega, \cF, \Pt)$. If the function $h$
is bounded, then there is a continuous embedding
\be
\label{emb}
L_2(\Omega, \Pt) \subset L_1(\Omega, \P).
\ee
Indeed, if $\xi\in L_2(\Omega, \Pt)$, then
\bee
\E \xi =\Et (Z_T \xi) \leq \sqrt{\Et Z_T^2} \sqrt{\Et \xi^2} \leq
C\sqrt{\Et \xi^2},
\bee
because
\bee
\ba{lll}
\Et Z_T^2 &=& \Et \left( \exp\left\{\int_0^T|h(X(t))|^2dt \right\}
 \exp\left\{ 2\int_0^T h^*(X(t))dY(t)-2\int_0^T|h(X(t))|^2dt \right\}
 \right)\\
&\leq& C \Et \exp\left\{ 2\int_0^T h^*(X(t))dY(t)-2\int_0^T|h(X(t))|^2dt
 \right\}
\leq C
\ea
\bee
where the last inequality follows from the property {\bf P2}
 of $\Pt$ and
Proposition 3.5.12 in \cite{KarShr}.

Next, consider the partial differential operators
\bee
\L g(x)=\frac{1}{2} \sum_{i,j=1}^d \left( (\sg(x)\sg^*(x))_{ij} +
(\rho(x)\rho^*(x))_{ij} \right) \frac{\pd^2 g(x)}{\pd x_i \pd x_j}+
\sum_{i=1}^d b_i(x) \frac{\pd g(x)}{\pd x_i};
\bee
\bee
\cM_l g(x) = h_l(x)g(x)+\sum_{i=1}^d \rho_{il}(x)
\frac{\pd g(x)}{\pd x_i},\
l=1, \ldots, r;
\bee
and their adjoints
\bee
\L^* g(x)\!=\!\frac{1}{2} \sum_{i,j=1}^d \frac{\pd^2}{\pd x_i \pd x_j}
\left( (\sg(x)\sg^*(x))_{ij}g(x) \!+ \!
(\rho(x)\rho^*(x))_{ij} g(x)\right) -
\sum_{i=1}^d \frac{\pd}{\pd x_i}\left( b_i(x)g(x) \right);
\bee
\bee
\cM_l^* g(x) = h_l(x)g(x)-\sum_{i=1}^d \frac{\pd}{\pd x_i}\left(
\rho_{il}(x)g(x) \right),\  l=1, \ldots, r.
\bee
Let $\H^a$ be the Sobolev space $\{f:(1+|w|^2)^{a/2}\hat{f}\in L_2(\R^d)\}$,
where $\hat{f}=\hat{f}(w)$ is the Fourier transform of $f$; $\H^0=L_2(\R^d)$
with the norm $\|\cdot\|_0$. The inner product in $L_2(\R^d)$ and the duality between
$\H^1$ and $\H^{-1}$ will be denoted by  $(\cdot, \cdot)_0.$
 Note that
the operators $\L, \L^* $ are bounded from $\H^1$ to $\H^{-1}$,
operators $\cM, \cM^*$ are bounded from $\H^1$ to $L_2(\R^d)$, and,
for every $g\in \H^1$,
\be
\label{eq:dissip}
2(\L^*g,g)_0+\sum_{l=1}^r \|\cM_l^*g\|_0^2 \leq C\|g\|_0^2.
\ee

The following result is well known \cite[Theorem 6.2.1]{Roz}.

\begin{prop}
\label{prop:egun}
In addition to Assumption R1 suppose that the initial density $p_0$
belongs
to the space $\H^1$. Then
there is a random field
$p=p(t,x),\ t \in [0,T],\ x \in \R^d,$ with the following properties:

1. $p\in L_2(\Omega\times (0,T), d\Pt\times dt; \H^1) \cap
L_2(\Omega, \Pt; \C([0,T], L_2(\R^d))). $

2. The function   $p(t,x)$ is a generalized solution
of  the stochastic partial differential equation
\be
\label{ze}
\ba{ll}
\ds dp(t,x)&\ds =\L^* p(t,x)dt+\sum_{l=1}^r\cM_l^*p(t,x) dY_l(t),\ 0<t\leq T,\
 x \in \R^d;\\
\ds p(0,x)&\ds =p_0(x).
\ea
\ee

3. The equality
\be
\label{uof}
\Et \left[ f(X(t))Z_t|\cF_t^Y \right] = \int_{\R^d} f(x) p(t,x) dx
\ee
holds for all bounded measurable functions $f$.
\end{prop}

The random field $p=p(t,x)$ is called the {\em unnormalized
filtering density}
(UFD) and the random variable
$\phi_t[f]=\Et \left[ f(X(t))Z_t|\cF_t^Y \right]$, the
{\em unnormalized
optimal filter}. Under Assumption R1, equation (\ref{ze}) is at least dissipative.
If the matrix $\sigma \sigma^*$ is uniformly positive
definite, then equation (\ref{ze}) is coercive rather than dissipative, and
it is enough to assume that $p_0 \in L_2(\R^d)$.

\section{Approximation of the optimal filter}
\setcounter{equation}{0}

Let $\{e_i,\ i\geq 1\}$ be an orthonormal basis in $L_2(\R^d)$ so that
every $e_i$ belongs to $ \H^1$. Fix a positive integer number $K$.
Define the matrices $\ds A^{K}=(A^{K}_{ij},\ i,j=1,\ldots,K)$
and $ \ds B^{K}_l=(B^{K}_{l,ij},\ i,j=1,\ldots,K;\ l=1, \ldots, r),$ by
\bee
A^{K}_{ij}=(\L^* e_{j}, e_{i})_0,\quad
B^{K}_{l,ij}=(\cM_l^* e_{j}, e_{i})_0.
\bee
Since $e_{i} \in \H^1$ for all $i$, the matrices  are well
defined.
The Galerkin approximation $p^{K}(t,x)$ of $p(t,x)$ is given by
\be
\label{eq:ga}
p^{K}(t,x)=\sum_{i=1}^K p^{K}_i(t) e_{i}(x),
\ee
where the vector $ p^{K}(t)=\{ p^{K}_{i}(t), \ i=1, \ldots, K\}$ is
the solution of the system of stochastic ordinary differential
equations
\be
\label{ga.sss2}
d p^{K}(t) =  A^{K} p^{K}(t) dt + \sumr B^{K}_lp^{K}(t) dY_l(t)
\ee
with the initial condition $p^{K}_{i}(0) = (p_0, e_{i})_0.$
Note that the matrices $  \ds B^{K}_l,\ l=1, \ldots, r,$ do not,
 in general,
commute with each other even if $\rho(x) \equiv 0$.

We next use Theorem \ref{Th.sp1} to
derive the Cameron-Martin version of the Wiener chaos
expansion of the solution of (\ref{ga.sss2}).

Let $0=t_0< t_1 \ldots <t_M=T$ be a uniform (for simplicity)
partition of the interval $[0,T]$
with step $\D$ and let $\{m_k(t),\ k\geq 1\}$ be an orthonormal basis
in $L_2([0,\D])$.  Denote by $J$ the set of all multi-indices
$\a=\{\a^l_k,\ l=1,\ldots, r,\ k\geq 1, \a^l_k=0,1,2,\ldots\}$ so that
$|\a|=\sum_{l,k} \a^l_k < \infty$.

Define random variables
\be
\label{eq:ch1}
\xi^i_{k,l}=\int_{t_{i-1}}^{t_i}
m_k(s-t_{i-1})dY_l(s),
\ee
and then, for $\a \in J$,
\be
\label{eq:ch2}
\xi^i_{\a} = \frac{1}{\sqrt{\a!}}
\prod_{k,l} H_{\a^l_k} (\xi^i_{k,l}),
\ee
where $ \ds H_n(t)=(-1)^n e^{t^2/2} \frac{d^n}{dt^n} e^{-t^2/2}$.

The following result is a direct consequence of
Theorem \ref{Th.sp1}.

\begin{theorem}
\label{th:wce}
 For every $i=1,\ldots, M$, the solution of (\ref{ga.sss2}) can be
 written in $L_2(\Om; \R^K)$ as
\be
\label{eq:wce}
p^{K}(t_i)=\sum_{\a \in J} \frac 1{\sqrt{\a !}}
\f_{\a}^{K}(\D;p^{K}(t_{i-1} ))\xi_{\a}^i,\ \ i=1,\ldots,M,
\ee
where, for $s\in (0,\D]$ and $\zeta\in \R^K$, the functions
 $\f_{\a}^{K}(s;\zeta)$ are the solutions of
\be
\label{eq:wcecoef}
\begin{array}{lll}
\ds \frac{\pd \f_{\a}^{K}(s;\zeta)}{\pd s}\!\!\!\!
 &=\ds A^{K}\f_{\a}^{K}(s;\zeta)+\sum_{k,l}
 \a_k^l m_k(s) B^{K}_l \f_{\a(k,l)}^{K}(s;\zeta), \  0<s \leq \D, \\
\ds \f_{\a}^{K}(0;\zeta)\!\!\!\! &= \ds \zeta 1_{ \{ | \a |=0 \} },
\end{array}
\ee
and
 $\a(i,j)$ stands for the  multi-index
 $\tilde{\a} = (\tilde{\a}_k^l)_{1 \leq l \leq r,\ k \geq 1}$ with
\be
\label{al1}
\tilde{\a}_k^l=
\begin{cases}
\a_k^l & {\rm \ if} \  k \not= i {\rm \ or}\  l \not= j {\rm \ or\  both} \cr
\max(0,\a_i^j-1) &{\rm \ if\ } k=i {\rm \ and\ } l=j.
\end{cases}
\ee
\end{theorem}

For fixed positive integers $N$ and $n$ define the set $J^n_N$
as the collection of multi-indices $\a$ from $J$ such that
$|\a|\leq N$ and
$\a^l_k=0$ if $k>n$. The approximation
$p^{K,n}_N(t_i)$ of $p^{K}(t_i)$ is defined by
\be
\label{eq:wce_appr}
p^{K,n}_N(t_0)=p^{K}(0),\quad
p^{K,n}_N(t_i)=\sum_{\a \in J_N^n}  \frac 1{\sqrt{\a !}}
\f_{\a}^{K}(\D;p^{K,n}_N(t_{i-1}))\xi_{\a}^i,\ \ i=1,\ldots,M.
\ee

\ni Note the $p^{K,n}_N(t_i)$ is a vector in $\R^K$.
Let $\bU=\{\bu^{j}, \ j=1, \ldots, K\}$ be a
 basis in $\R^{K}$.  The vector $p_N^{K,n}(t_i)$
can then be written as
\bee
 p_N^{K,n}(t_i) = \sum_{j=1}^K p_{N,j}^{K,n}(t_i;\bU)\bu^j,
\bee
and by the recursive definition of $p_N^{K,n}(t_i)$,
\bee
\ba{lll}
\ds  p_N^{K,n}(t_{i+1}) &=& \ds \sum_{\a \in J^n_N}
\f_{\a}^{K}(\D;p_N^{K,n}(t_{i})) \xi_{\a}^i\\
&=&\ds  \sum_{\a \in J^n_N}\sum_{j=1}^K
\f_{\a}^{K}(\D;\bu^j)p_{N,j}^{K,n}(t_i;\bU)  \xi_{\a}^i.
\ea
\bee
Once again, $\varphi^k_{\a}(\Delta,\bu^i)$ is a vector in $\R^K$,
so  we write
\bee
\f_{\a}^K(\D, \bu^j)=\sum_{k=1}^K q^{K,\a}_{jk}(\bU) \bu^k,
\bee
 and conclude that
\be
\label{papr.sss2}
p_{N,j}^{K,n}(t_{i+1};\bU) =   \sum_{\a \in J^n_N}\sum_{k=1}^K
q^{K,\a}_{jk}(\bU) p_{N,k}^{K,n}(t_{i};\bU) \xi^i_{\a}.
\ee
Then
\be
\label{eq:ufda}
p_{N}^{K,n}(t_i,x)=
\sum_{j,k=1}^Kp_{N,j}^{K,n}(t_{i+1};\bU)\bu^j_{k} e_{k}(x)
\ee
is an approximation of the unnormalized filtering density.

Suppose that the basis functions $e_k$ and the function $f$ are such
that
\be
\label{eq:fcoef}
f_k=\int_{\R^d} f(x)e_k(x)dx
\ee
is defined for every $k=1,\ldots, K$. It follows from (\ref{eq:ufda})
that
\be
\label{eq:uofa}
\tilde{\phi}_i[f]=\sum_{j,k=1}^Kp_{N,j}^{K,n}
(t_{i+1};\bU)\bu^j_{k} f_{k}
\ee
is an approximation of the unnormlized optimal filter.

The following is a possible algorithm for computing approximations of the
unnormlized filtering density and optimal filter using (\ref{eq:ufda}) and
(\ref{eq:uofa}).

{\it
1. $\underline{Preliminary \ computations}$ (before the observations
are available):
\begin{enumerate}
\item Choose suitable basis functions $\{e_k, k=1, \ldots, K\}$
in $L_2(\R^d)$,
$\{m_i, i=1, \ldots, n\} $ in $L_2([0,\Delta])$, and a {\bf standard unit
basis} $\{ \bu^j, j=1, \ldots K\}$ in $\R^K$, that is, $\bu_i^i=1$, $\bu_i^j=0$
otherwise.
\item
 for $\a\in J_N^n$ and $j,k=1,\ldots, K$ compute
 $$
 q^{K,\a}_{jk}= \f_{\a,j}^{K}(\D;{\bf u}^{k}) {\rm\  (using\
 (\ref{eq:wcecoef}))},
 \  f_{k}=\int_{\R^d}f(x)e_k(x)dx, \ p_{N,k}^{K,n}(t_0)=\int_{\R^d} p_0(x)e_k(x)dx;
$$
\end{enumerate}

2. $\underline{Real-time\  computations, \ i-th \ step}$
 (as the observations become available):  compute
 $\xi^i_{\a}$, $\a \in J^n_N$ (according to (\ref{eq:ch1})
 and (\ref{eq:ch2}));
 \bee
 Q_{jk}^{K}(\xi^i)= \sum_{\a \in J_N^n}q^{K,\a}_{jk}\xi_{\a}^i;
\bee
\be
\label{psi.sss2}
p_{N,j}^{K,n}(t_{i}) = \sum_{k=1}^K
 Q_{jk}^{K}(\xi^i)p_{N,k}^{K,n}(t_{i-1}), \quad  j=1, \ldots, K;
\ee

then, if necessary, compute
\be
\label{pb.sss2}
p^{K,n}_N(t_i,x)=\sum_{j=1}^Kp_{N,j}^{K,n}(t_{i})  e_{j}(x),
\ee
\be
\label{fb.sss2}
\tilde{\phi}_{t_i}[f]=
\sum_{j=1}^K  f_{j}p_{N,j}^{K,n}(t_{i}),
\ee
and
\be
\label{nof.sss2}
\tilde{f}_{t_i}=
\frac{\tilde{\phi}_{t_i}[f]}{\tilde{\phi}_{t_i}[1]}.
\ee
}

\br
\label{sss1.rm2} The main advantage of the above algorithm
  as compared to most other
schemes for solving the Zakai equation is that the time consuming
computations, including solving  partial differential equations and
computing
integrals, are performed in advance, while the real-time  part is
relatively
simple even when the dimension $d$ of the state process is large.
Here are some other features of the algorithm:
\begin{enumerate}

\item The overall amount of preliminary
 computations does not depend on the number of the on-line time steps;
\item Formulas (\ref{fb.sss2}) and (\ref{nof.sss2})
can be used to compute an approximation to $\hat{f}_{t_i},$
for example,  conditional moments, without the time consuming computations of
$p^{\k,n}_N(t_i,x)$ and the related integrals;
\item Only the  coefficients $p^{K,n}_{N,j}(t_i)$  must be computed
at every time step while the approximate filter $\tilde{f}_{t_i}$ and/or
UFD $p^{K,n}_N(t_i,x)$ can be computed as needed, for example, at the final time
moment.
\item The real-time  part of the algorithm can be easily parallelized.
\item Even though the coefficients $q^{K,\a}_{jk}$ are computed
according to (\ref{eq:wcecoef}), their values can be further adjusted by
simulating the state and observation processes and computing the corresponding
filter estimates.
\item If  $n=1$, then each
$\xi_{\a}^i$ depends only on the increments $Y_l(t_{i})-Y_l(t_{i-1})$ of the
observation process. For $n>1$  and $k>1$,  the integral
$\ds \int_{t_{i-1}}^{t_{i}} m_k(s-t_{i-1})dY_l(s)$ can be reduced to
a usual Riemann integral and then approximated
by the trapezoidal  rule.
\end{enumerate}
\er

\section{Rate of convergence}

To study the convergence of the algorithm, it is necessary to
specify the bases $\{e_k, k\geq 1\}$ on $\R^d$ and $\{m_i, i\geq 1\}$
on $[0,\D]$.

Let  $\{e_{k}, k\geq 1\}$  be the Hermite basis in $L_2(\R^d)$.
The basis can be described as follows. Denote by $\Gm$ the set of ordered
$d$-tuples $\gm=(\gm_1, \ldots, \gm_d)$ with $\gm_j=0,1,2,\ldots$.
For $\gm \in \Gm$ define
\bee
\cH_{\gm}(x)=\prod_{j=1}^d \cH_{\gm_j}(x_j),
\bee
where
\bee
\cH_k(t)=\frac{(-1)^n}{\sqrt{2^n\pi^{1/2}n!}} e^{t^2}\frac{d^n}{dt^n}
e^{-t^2}
\bee
With this definition,
$\cH_{\gm}$ is the eigenfunction of the self-adjoint operator
$\Ld=-\nabla^2+(1+|x|^2)$:
\bee
\Ld \cH_{\gm} = \ld_{\gm}e_{\gm},
\bee
where $\nabla^2$ is the Laplace operator and $\ld_{\gm}=(2|\gm|+d+1)$.

To define an ordering of the set $\Gm$, we
 define $|\gm|=\sum_{j=1}^d\gm_j$ and then say that $\gm<\tau$ if
$|\gm|<|\tau|$ or if $|\gm|=|\tau|$ and $\gm<\tau$ under the
lexicographic
ordering, that is, $\gm_{i_0} < \tau_{i_0}$, where $i_0$ is the
first index
for which $\gm_i\not= \tau_i$. The basis  $\{e_{k}\}_{k\geq 1}$
is then the set $\{ \cH_{\gm}(x), \ \gm \in \Gm\}$ together with the
above ordering of the set $\Gm$ so that $\Ld e_k =\ld_k e_k$ and
$\ld_k \asymp k^{1/d}.$

Next, we define an orthonormal basis $\{m_k\}$ in $L_2([0,\D])$
by
\bee
 m_1(s)=\frac 1{\sqrt{\D}}; \quad m_k(s)=\sqrt{\frac{2}{\D} }
 \cos \left( \frac{\pi (k-1) s}{\D} \right), \, k>1; \ 0\leq s \leq \D.
\bee

\begin{definition}
\label{regfm.d}
The filtering model (\ref{fmod.sss2}) is called $\nu$-regular for
some positive integer $\nu$ if the functions $\sigma$ and $\rho$ belong to
$\C^{2\nu+3}_b$, the functions  $b$ and $h$ belong to $\C^{2\nu+2}_b$,
and $\Ld^{\nu} p_0 \in \H^1$.
\end{definition}

\bt
\label{mth.sss2}
If the filtering model (\ref{fmod.sss2}) is $\nu$-regular, in the sense of
Definition \ref{regfm.d}, for some $\nu >d+1$ and
\bee
C_{\rho}=\max_{i,l} \sup_{x \in \R^d} |\rho_{il}(x)|^2,
\bee
then
\be
\label{ufde.sss2}
\ba{l}
\ds \max_{0 \leq i \leq M}
\Et\|p(t_i,\cd)-p^{K,n}_N(t_i,\cd)\|_0^2 \leq
\frac{C(\nu,T)}{K^{2(\nu-d-1)/d}} \\ +
\ds \left( C \frac{(1+C_{\rho} K^{1/d})\D+(K^{2/d}+
C_{\rho}K^{3/d})\D^2}{n}+
\frac{(C(1+C_{\rho}K^{1/d}))^{N+1} \D^N}{(N+1)!}
\right) e^{C(1+C_{\rho}K^{1/d})T}.
\ea
\ee
The number $C(\nu,T)$ depends on $\nu, T$, and the parameters of the model
(coefficients of the equations (\ref{fmod.sss2})). The number $C$ depends only
on the parameters of the model.

If, in addition, $(1+|x|^2)^{-w}f \in L_2(\R^d)$ for some
$w\geq 0$ so that
$\nu>d+1+w$ and \\
$\Ld^{\nu}((1+|x|^2)^{w}p_0)\in \H^1$, then
\be
\label{uofe.sss2}
\ba{l}
\ds \max_{0 \leq i \leq M}
\Et|\phi_{t_i}[f]-\tilde{\phi}_{t_i}[f]|^2 \leq
\frac{C(\nu,T,w)C_f}{K^{2(\nu-w-d-1)/d}} \\+
\ds C_f\left( C \frac{(1+C_{\rho} K^{1/d})\D+(K^{2/d}+
C_{\rho}K^{3/d})\D^2}{n}+
\frac{(C(1+C_{\rho}K^{1/d}))^{N+1} \D^N}{(N+1)!} \right)
e^{C(1+C_{\rho}K^{1/d})T}.
\ea
\ee
The number $C(\nu,T,w)$ depends on $\nu, T,w$, and the parameters of the model;
 the number $C$ depends only on $w$ and the parameters of the model;
 $C_f=\int_{\R^d} (1+|x|^2)^{-2w} |f(x)|^2dx.$
\et

{\bf Proof.} By Theorem \ref{H2.th},
\be
\label{eq:pr_ga}
\Et\|p(t_i,\cdot)-p^K(t_i, \cdot)\|^2_0 \leq
\frac{C(\nu,T)}{K^{2(\nu-d-1)/d}}.
\ee
Indeed, by Theorem 4.3.2 in \cite{Roz},
$\sup_{0<t<T}\Et\|\Ld^{\nu}p(t,\cdot)\|_0^2
\leq e^{CT}\|\Ld^{\nu}p_0\|_0^2$, where $C$ depends only on
$\nu$ and the parameters of the model. Also, in the notations
of Theorem \ref{H2.th}, $\theta=1/d$, $q=1/(2d)$, and
$\theta_1=(\nu-1)/d$.

 To simplify the further presentation, set  $\k=K^{1/d}$ and
 define $C_{\k}=1+C_{\rho}\k$.
Then, to prove (\ref{ufde.sss2}), it remains to show that
\bee
\Et|p^{K}(t_i)-p^{K,n}_N(t_i)|^2 \leq
\left( C \frac{C_{\k}\D+\k^2C_{\k}\D^2}{n}+
\frac{(CC_{\k})^{N+1} \D^N}{(N+1)!} \right) e^{CC_{\k}T},
\bee
and by Theorem \ref{utrn.cor} this inequality holds if, for
every vector
$\zeta \in \R^{K}$,
\be
\label{A.sss2}
|A^{K}\zeta|^2 \leq C\k^2 |\zeta|^2,\
|B_l^{K} \zeta|^2 \leq CC_{\k}|\zeta|^2,\
|e^{tA^K} \zeta|^2 \leq e^{Ct}|\zeta|^2.
\ee
Because of the multi-step approximation, we, as usual, loose one power of
$\Delta$ in (\ref{eq:trNn}).
Inequalities (\ref{A.sss2}) are verified  by direct calculations
using  that the operators $\Ld^{-1}\L$ and $\Ld^{-1/2}\cM_l$
are bounded in $L_2(\R^d)$.

To prove (\ref{uofe.sss2}), let $\beta(x)=\sqrt{1+|x|^2}$ and,
for $w \in \R$, define the space $L_{2,w}(\R^d)=
\{ f: f\beta^w \in L_2(\R^d) \}.$ Clearly, $L_{2,w}(\R^d)$ is a
Hilbert space with inner product $(f,g)_{0,w}=(f\beta^w, g\beta^w)_0$
 and norm $\|f\|_{0,w}^2=(f,f)_{0,w}$. Then
 \be
 \label{eq:pr202}
 |\phi_{t_i}[f]-\tilde{\phi}_{t_i}[f]|^2\leq C_f \|p(t_i,\cdot)-
 p^{K,n}_N(t_i,\cdot)\|_{0,2w}^2.
 \ee
 Using the calculus of pseudo-differential operators \cite[Chapter 4]{Shu},
 we conclude that, for every $g \in L_{2,2w}(\R^d)$,
 \be
 \label{eq:pr200}
 \|g\|_{0,2w}=\|\beta^{2w}g\|_0\leq C\|\Ld^{-w}\beta^{2w}\Ld^{w}g\|_0 \leq
 C\|\Ld^{w}g\|_0.
 \ee
 Therefore, by Theorem \ref{H2.th} and Theorem 4.3.2 in \cite{Roz},
\be
\label{eq:pr_ga2}
\Et\|p(t_i,\cdot)-p^K(t_i, \cdot)\|^2_{0,2w} \leq
\frac{C(\nu,T)}{K^{2(\nu-w-d-1)/d}}.
\ee
Next,  (\ref{eq:pr200}) implies
\bee
\Et \|p^{K}(t_i,\cd) - p^{K,n}_N(t_i, \cd)\|^2_{0,2w} \leq
C\sum_{k=1}^K \ld_{k}^{2w}
\Et| p^{K}_{k}(t_i)- p^{K,n}_{N,k}(t_i)|^2.
\bee
Define the diagonal matrix
$\ds \hat{\Ld}=(\hat{\Ld}_{ij})_{i,j=1, \ldots, K}$
by $\ds \hat{\Ld}_{ii}=\ld_{i}^{w}$. Then define the matrices
\bee
\hat{A}^K = \hat{\Ld}A^{K} \hat{\Ld}^{-1},\
\hat{B}_l^K = \hat{\Ld}B^{K}_l \hat{\Ld}^{-1},
\bee
and the vectors $\ds \hat{p}(t)=\hat{\Ld}p^{K}(t),\ \hat{p}^n_N(t_i)=\hat{\Ld}
p^{K,n}_N(t_i)$. With these definitions,
the vector
 $\hat{p}^K(t)$ is the solution of
\bee
d \hat{p}^K(t) =  \hat{A}^K\hat{p}^K(t) dt + \sumr \hat{B}_l^K\hat{p}(t) dY_l(t)
\bee
with the initial condition
$\ds\hat{p}_{k}^K(0)=\ld^{w}_{k}(p_0,e_{k})_0,$
 the vector $\hat{p}^{K,n}_N(t)$ satisfies
\bee
\hat{p}^{K,n}_N(t_0)=\hat{p}^K(0),\
\hat{p}^{K,n}_N(t_i)=\sum_{\a \in J_N^n}  \frac 1{\sqrt{\a !}}
\hat{\f}_{\a}^{K}(\D;\hat{p}^{K,n}_N(t_{i-1}))\xi_{\a}^i,\ \ i=1,\ldots,M,
\bee
and
\be
\label{eq:pr205}
\Et \|p^{K}(t_i,\cd) - p^{K,n}_N(t_i, \cd)\|^2_{0,w}\leq
C\Et|\hat{p}^K(t_i)-\hat{p}^{K,n}_N(t_i)|^2.
\ee
The functions $\hat{\varphi}^K_{\alpha}$ satisfy the equations (\ref{eq:wcecoef})
with $\hat{A}^K$ and $\hat{B}_l^K$ instead of $A^K$ and $B^K$.

Direct computations show that, for all $\zeta \in \R^{K}$,
\be
|\hat{A}^{K}\zeta|^2 \leq C\k^2 |\zeta|^2,\  \  \
|\hat{\B}_l^{K} \zeta|^2 \leq CC_{\k}|\zeta|^2, \  \  \
|e^{t\hat{A}^{K}} \zeta|^2 \leq e^{Ct}|\zeta|^2,
\ee
with $C$ depending on $w$ and the parameters of the filtering model.
By Theorem \ref{utrn.cor} we then conclude that
\bee
\Et|\hat{p}^{K}(t_i)-\hat{p}^{K,n}_N(t_i)|^2 \leq
\left( C \frac{C_{\k}\D+\k^2C_{\k}\D^2}{n}+
\frac{(CC_{\k})^{N+1} \D^N}{(N+1)!} \right) e^{CC_{\k}T}.
\bee
Together with (\ref{eq:pr202}), (\ref{eq:pr200}), and (\ref{eq:pr205}),
the last inequality implies (\ref{uofe.sss2}).

Theorem \ref{mth.sss2} is proved.
\endp

\section{Proof of Theorem \protect{\ref{utrn.cor}}}
\setcounter{equation}{0}

The proof requires an explicit formula for the solution of (\ref{S}).
We begin with some auxiliary constructions.

Every multi-index $\a$ with $|\a|=k$ can be identified with the set $K_{\a}=
\{ (i_1^{\a},q_1^{\a}), \ldots, (i_k^{\a},q_k^{\a}) \}$ so that
$i_1^{\a} \leq i_2^{\a}\leq \ldots \leq i_k^{\a}$ and if
$i_j^{\a}=i_{j+1}^{\a}$, then $q_j^{\a} \leq q_{j+1}^{\a}$. The
first pair $(i_1^{\a},q_1^{\a})$ in $K_{\a}$ is the position numbers of the
first nonzero element of $\a$. The second pair is the same as the first
if the first nonzero element of $\a$ is greater than one; otherwise,
the second pair is the position numbers of the second nonzero element
 of $\a$ and so on. As a result, if $\a^q_j>0$,
 then  exactly $\a_j^q$ pairs in $K_{\a}$
are $(j,q)$. The set $K_{\a}$ will be referred to as {\bf the characteristic
set} of the  multi-index $\a$. For example, if $r=2$ and
\bee \index{Characteristic set of a multi-index}
\a=\left(
\ba{llllllll}
0 &1 &0 &2 &3 &0 &0 &\ldots \\
1 &2 &0 &0 &0 &1 &0 &\ldots
\ea
\right),
\bee
then the nonzero elements are $\a_1^2=\a_2^1=\a_1^6=1,\ \a_2^2=\a_4^1=2,\
 \a_5^1=3,$ and the \\ characteristic set is
$K_{\a}\!=\!\{ (1,2),\,(2,1),\,(2,2),\,(2,2),\,(4,1),\,(4,1),\,
(5,1),\,(5,1),\,(5,1),\,(6,2) \}$.
In the future, when there is no danger of confusion, the superscript
$\a$ in $i$ and $q$ will be omitted so that  $(i_j,q_j)$ will be written  instead of
$(i_j^{\a},q_j^{\a})$.

Let $\PP^k$ be the permutation group of the set $ \{ 1, \ldots, k \} $.
For a given $\a \in J$ with $| \a|=k$ and the characteristic set
 $\{(i_1,q_1),\ldots, (i_k,q_k)\}$
define
\bee
E_{\a}(s^k;l^k):= \sum_{\sigma \in \PP^k}
 m_{i_1}(s_{\sigma(1)})1_{ \{ l_{\sigma(1)}=q_1 \} } \cdots
m_{i_k}(s_{\sigma(k)})1_{ \{ l_{\sigma(k)}=q_k \} }.
\bee

The following notations are introduced to simplify the further presentation:
\begin{itemize}
\item $s^k$, the  ordered set $(s_1, \ldots, s_k)$; $ds^k:=ds_1\ldots ds_k$;
\item $l^k$, the  ordered set $(l_1, \ldots, l_k)$;
\item $\F_t=e^{At}$;
\item  $F(t;s^k;l^k;g) := \F_{t-s_k}B_{l_k}
\F_{s_k-s_{k-1}} \ldots B_{l_1} \F_{s_1-T_0}g, \ k \geq 1;$
\item $\ds  \int_{T_0}^{(k,t)} (\cdots) ds^k := \int_{T_0}^t \int_{T_0}^{s_k}
\ldots \int_{T_0}^{s_2} (\cdots) ds_1\ldots ds_k$;
\item  $\ds \sum\limits_{l^k} :=\sum\limits_{l_1,\ldots, l_k=1}^r$.
\end{itemize}
Note that
\be
\label{eq:pr_tech}
|F(t;s^k;l^k;g)|^2 \leq C_2^k e^{C_1(t-T_0)}|g|^2,\
\int_{T_0}^{(k,t)} ds^k=\frac{(t-T_0)^k}{k!}, \
\sum\limits_{l^k}1=r^k.
\ee

\ni \begin{prop}
\label{S.th}  If  $\a \in J$ is a multi-index  with $|\a|=k$ and the characteristic
set $\{ (i_1,q_1),\ldots, (i_k,q_k) \}$, then, for
$t \in [T_0,t^*]$, the  corresponding solution
 $ \f_{\a}(t;T_0;U_0) $ of (\ref{S}) is given by
\be
\label{A15}
\begin{array}{l}
\displaystyle \f_{\a}(t;T_0;U_0)=\\
\ds \sum_{\sigma \in \PP^k}\sum_{l^k} \int\limits_{T_0}^{(k,t)}
 F^{k}(t;s^k;l^k;U_0) m_{i_{\sigma(k)}}(s_k)1_{ \{ l_k=q_{\sigma(k)} \} }
\cdots m_{i_{\sigma(1)}}(s_1)1_{ \{ l_1=q_{\sigma(1)} \} } ds^k,\ k\!>\!1;\\
\displaystyle \f_{\a}(t;T_0;U_0)=
\int_{T_0}^t \F_{t-s_1}B_{q_1}\F_{s_1-T_0}U_0m_{i_1}(s_1)ds_1, \
k=1;\\
\displaystyle \f_{\a}(t;T_0;U_0)=\F_{t-T_0}U_0,\ k=0,
\end{array}
\ee
and
\be
\label{A20}
\sum_{|\a|=k} \frac{|\f_{\a}(t;T_0;U_0)|^2}{\a!} =
\sum_{l^k}\int_{T_0}^{(k,t)} |F(t;s^k;l^k;U_0)|^2 ds^k.
\ee
 \end{prop}

\ni \pr To simplify the notations,
 the arguments  $T_0$ and $U_0$  will be omitted wherever possible.
 Representation (\ref{A15}) is obviously true for $|\a|=0$.
Then the general case $|\a|\geq 1$ follows by induction from
 the variation of parameters formula.

To prove (\ref{A20}), first of all note that
\bee
\ba{l}
\ds \sum_{\sigma \in \PP^k} m_{i_{\sigma(k)}}(s_k)1_{ \{ l_k=q_{\sigma(k)} \} } \cdots
m_{i_{\sigma(1)}}(s_1)1_{ \{ l_1=q_{\sigma(1)} \} } \\ =
\ds \sum_{\sigma \in \PP^k} m_{i_k}(s_{\sigma(k)})1_{ \{ l_{\sigma(k)}=q_k \} } \cdots
m_{i_1}(s_{\sigma(1)})1_{ \{ l_{\sigma(1)}=q_1 \} }.
\ea
\bee
 Indeed, every term on the left corresponding to a given $\sigma_0 \in \PP^k$
coincides with
 the term on the right corresponding to $\sigma^{-1}_0 \in \PP^k$.

Then  (\ref{A15}) can be written as
$\f_{\a}(t)=\sum_{l^k}\int_{T_0}^{(k,t)}F(t;s^k;l^k)E_{\a}(s^k;l^k)ds^k.$
 Using the notation
\bee
G(t;s^k;l^k) := \sum_{\sg \in \PP^k}
\F_{t-s_{\sg(k)}}B_{l_{\sg(k)}} \ldots
\F_{s_{\sg(2)}-s_{\sg(1)}}B_{l_{\sg(1)}}\F_{s_{\sg(1)}-T_0}g
1_{s_{\sg(1)}< \ldots < s_{\sg(k)}<t},
\bee
\ni it   can  be rewritten  as
\be
\label{A25}
\f_{\a}(t)= \frac{1}{k!}
\sum_{l^{k}}\int\limits_{[T_0,t^*]^k} G(t;s^k;l^k) E_{\a}(s^k;l^k)ds^k.
\ee

Since for every $t \in [T_0,t^*]$ the function   $G(t;s^k;l^k)$ is  symmetric,
$$
G(t;s^k;l^k)=\sum\limits_{|\beta|=k}\frac{c_{\beta}(t)
E_{\beta}(s^k;l^k)}{\sqrt{\beta! k!}}
$$
with some vector coefficients $c_{\beta}(t)$. This and  (\ref{A25}) imply
$ |\f_{\a}(t)|^2/\a!=|c_{\a}|^2/k!$ and so
\bee
\ba{c}
\ds\sum_{|\a|=k} \frac{|\f_{\a}(t)|^2}{\a!} =
\frac{1}{k!}\sum_{|\a|=k}|c_{\a}(t)|^2 =
 \frac{1}{k!} \int\limits_{[T_0,t^*]^k} |G(t;s^k;l^k)|^2ds^k\\ =
\ds \frac{1}{k!}\sum_{l^k} \int\limits_{[T_0,t^*]^k} \Big{|}
 \! \sum_{\sg \in \PP^k}
\!\F_{t-s_{\sg(k)}}B_{l_{\sg(k)}} \ldots
 \F_{s_{\sg(2)}-s_{\sg(1)}}B_{l_{\sg(1)}}\F_{s_{\sg(1)}-T_0}g
1_{s_{\sg(1)}< \ldots < s_{\sg(k)}<t}
 \Big{|}^2 ds^k\\ =
\ds\sum_{l^k} \int_{T_0}^{(k,t)} |F(t;s^k;l^k)|^2 ds^k,
\ea
\bee
 which proves (\ref{A20}). Proposition \ref{S.th} is proved.
\endp

\ni  We continue by considering  the truncation only of the
length of $\a$.
Define $J_N=\{\a \in J:|\a| \leq N\}$ and
\be
\label{utrN.1}
U^N(s;T_0;U_0)=\sum_{\a\in J_N}\frac 1{\sqrt{\a
!}}\f_{\a}(s;T_0;U_0)\xi_{\a}(W_{T_0,t^*}).
\ee
Note that the summation in (\ref{utrN.1}) is still infinite.

\ni \begin{prop}
\label{utrN.th}
 In the notations of Theorem \ref{utrn.cor},
\be
\label{utrN.2}
\sup_{s \in [T_0,t^*]} \E |U(s;T_0;U_0)-U^N(s;T_0;U_0)|^2 \leq
\frac{[C_2r(t^*-T_0)]^{N+1}}{(N+1)!}e^{\Cb(t^*-T_0)}\;\E|U_0|^2.
\ee
\end{prop}

\ni \pr  To simplify the presentation, the arguments $T_0$ and
$U_0$ will be omitted wherever possible.

 By Theorem \ref{S.th},
\be
\label{pf:5}
\sum_{ | \a | = k} \frac{|\f_{\a}(s)|^2}{\a!} =
\sum_{l^k} \int_{T_0}^{(k,s)} |F(s;s^k;l^k)|^2 ds^k.
\ee
Since the random variables
$\xi_{\a}(W_{T_0,t})$ are uncorrelated and are independent of $U_0$,
formulas (\ref{eq:wce}) and (\ref{utrN.1}) imply
$\E | U(s) - U^N(s) |^2 =\sum_{k >N}\sum_{ | \a | = k}
\frac{\E|\f_{\a}(s)|^2}{\a!}.$
By  (\ref{eq:pr_tech}),
\bee
\ba{c}
\ds \sum_{k >N}\sum_{ | \a | = k} \frac{\E|\f_{\a}(s)|^2}{\a!} \leq
e^{C_1(s-T_0)}\;\E|U_0|^2 \sum_{k > N}
\frac{(C_2r(s-T_0))^k}{k!} \\ \leq
\ds \frac{(C_2r(t^*-T_0))^{N+1}}{(N+1)!}e^{\Cb(t^*-T_0)}\;\E|U_0|^2,
\ea
\bee
which completes the proof of Proposition \ref{utrN.th}.
\endp

Now we truncate the sum in (\ref{utrN.1}) even more by restricting $\alpha$ to the
set $J^n_N$.
\ni \begin{prop}
\label{prop:nN}
In the notations of Theorem \ref{utrn.cor} and Proposition \ref{utrN.th},
\be
\label{mer}
\ba{lll}
\ds \E|U^N(t^*;T_0;U_0)-U^n_N(t^*;T_0;U_0)|^2 &\leq&
\ds 2C_2r\;e^{\Cb(t^*-T_0)}\Big(
 \epsilon(B) \frac{(t^*-T_0)^2}{n}\\
&+&\ds C_0\left(1+(t^*-T_0)C_2r\right)\frac{(t^*-T_0)^3}{n}\Big)\E|U_0|^2.
\ea
\ee
\end{prop}

\ni \pr To simplify the presentation, the arguments $T_0$ and
$U_0$ will be omitted wherever possible.

If $\a$ is a multi-index with $|\a|=k$ and the characteristic set
$\{ (i_1^{\a},q_1^{\a}) \,\ldots, (i_k^{\a},q_k^{\a}) \}, $
 then $i_{k}^{\a}=d(\a)$, the order
of $\a$, and so the set $J_N^n$ can  be described as
$\{ \a \in J: |\a| \leq N; \ i_{|\a|}^{\a} \leq n \}$. Since
the random variables $\xi_{\a}$ are uncorrelated and are independent of $U_0$,
\bee
\E | U_N^n(t^*)-U^N(t^*)|^2 =\sum_{b=n+1}^{\infty} \sum_{k =1}^N
 \sum_{|\a|=k;i_k^{\a}=b} \frac{\E|\f_{\a}(t^*)|^2}{\a!}.
\bee

The problem is thus to estimate  $\ds \sum\limits_{b=n+1}^{\infty}
\sum\limits_{k = 1}^N\sum\limits_{|\a|=k;i_k^{\a}
 = b}\frac{|\f_{\a}(t^*)|^2}{\a!}.$

By Theorem \ref{S.th}  the corresponding solution $\f_{\a}$ of
(\ref{S}) can be written as
\be
\label{pf:20}
\f_{\a}(t^*)=\sum_{l^k}\int_{T_0}^{(k,t^*)}F(t^*;s^k;l^k)E_{\a}(s^k,l^k)ds^k.
\ee
According to  (\ref{al}), the characteristic set of $\a(i_k,q_k)$ is
 $\{ (i_1,q_1), \,\ldots, (i_{k-1},q_{k-1}) \}$;
\\ therefore,
it is possible to write
\bee
E_{\a}(s^k)=\sum_{j=1}^k m_{i_k}(s_j) 1_{ \{ l_j=q_k \} }
E_{\a(i_k,q_k)}(s^k_j;l^k_j),
\bee
 where $s^k_j$ (resp. $l^k_j$) denotes the same set $(s_1, \ldots, s_k)$
(resp. $(l_1, \,\ldots, l_k)$) with
omitted $s_j$ (resp. $l_j$); for example, $s^k_1=(s_2, \ldots, s_k).$

As a result, after changing  the order of integration
in the multiple integral, equality (\ref{pf:20}) can be rewritten  as
\be
\label{pf:25}
\f_{\a}(t^*)\! =\! \sum_{j=1}^k \sum_{l_j^k}
 \int\limits_{T_0}^{(k-1,t^*)}\!\! \Big(
\int_{s_{j-1}}^{s_{j+1}}\!\! F(t^*;s^k;l^k)m_{i_k}(s_j)1_{\{l_j=q_k\}}ds_j\Big)
E_{\a(i_k,q_k)}(s_j^k;l^k_j)ds_j^k,
\ee
\ni where $s_0:=T_0; \ s_{k+1}:=t^*$.

Denote
\bee
M_k(s) := \frac{\sqrt{2(t^*-T_0)}}{\pi(k-1)} \sin
\Big(\frac{\pi(k-1)(s-T_0)}{(t^*-T_0)} \Big);
\ k>1, \ T_0 \leq s \leq t^*,
\bee
\ni and  $\ds F_j := \frac{\pd F(t^*;s^k;l^k)}{\pd s_j}$. Then, as long as
$i_k=b > 1$,  integration by parts in the inner integral
on the right hand side  of (\ref{pf:25}) yields:
\bee
\ba{l}
\ds \int_{s_{j-1}}^{s_{j+1}} F(t^*;s^k;l^k)m_{b}(s_j)ds_j \\ =
\ds F(t^*;s^k;l^k)M_{b}(s_j) \Big|^{s_j=s_{j+1}}_{s_j=s_{j-1}} -
\int_{s_{j-1}}^{s_{j+1}} F_j(t^*;s^k;l^k)M_{b}(s_j)ds_j.
\ea
\bee

For each $j$, let us rename the remaining variables $s_j^k$ in
(\ref{pf:25}) as follows:
$t_i:=s_i, i \leq j-1; \ t_i:=s_{i+1}, i>j-1$, or, symbolically,
$t^{k-1}:=s_j^k$. We will set $t_0:=T_0$, $t_k:=t^*$ and
denote by $t^{k-1,j}, j=1, \ldots,k-1,$ the set
 $t^{k-1}$ in which $t_j$ is repeated twice (e.g.
$t^{k-1,1}=(t_1,t_1, \ldots,t_{k-1}), {\rm etc.})$; also
 $t^{k-1,0}:=(t_0,t_1,t_2,\ldots,t_{k-1})$,
 $t^{k-1,k}:=(t_1,\ldots,t_{k-1},t_k)$.

The similar changes will  also be made with the set $l^k$: for fixed $j$,
there are $k-1$ free indices $l_1 ,\ldots, l_{j-1}, l_{j+1} ,\ldots, l_k$ and
they are renamed just like $s^k$ to form the set $l^{k-1}$ (in this case,
the same symbols are used). Similarly, $l^{k-1,j}$ denotes the set
 $(l_1 ,\ldots, l_{j-1}, q_k, l_{j} ,\ldots, l_{k-1}).$ After these
transformations, $E_{\a(i_k,q_k)}(s_j^k;l^k_j)$ becomes
$E_{\a(i_k,q_k)}(t^{k-1};l^{k-1})$ - independent of $j$, and

\bee
\ba{l}
\ds F(t^*;s^k;l^k)1_{ \{ l_j=q_k \} }M_{b}(s_j)
 \Big|^{s_j=s_{j+1}}_{s_j=s_{j-1}} \\ =
 \ds F(t^*;t^{k-1,j};l^{k-1,j})M_{b}(t_j)-
 F(t^*;t^{k-1,j-1};l^{k-1,j})M_{b}(t_{j-1}),\
j=1, \ldots, k.
\ea
\bee
Therefore, if $d(\a)=b>1$ and $|\a|=k>0$, then
\bee
\ba{lll}
\ds \f_{\a}(t^*)&=&\ds \sum_{l^{k-1}} \int\limits_{T_0}^{(k-1,t^*)}
 \Big( f^{(1)}_b(t^*; t^{k-1}; l^{k-1})\\
&+&\ds f^{(2)}_b(t^*;t^{k-1}; l^{k-1})
 \Big)E_{\a(i_k,q_k)}(t^{k-1};l^{k-1})dt^{k-1},
\ea
\bee
where
\bee
\ba{lll}
  f^{(1)}_b(t^*; t^{k-1}; l^{k-1})&=&\ds \sum_{j=1}^k
 \Big(F(t^*;t^{k-1,j};l^{k-1,j})M_{b}(t_j) \\
 &-& \ds F(t^*;t^{k-1,j-1};l^{k-1,j})M_{b}(t_{j-1}) \Big)\quad {\rm if}\ k>1,
\ea
\bee

$f^{(1)}_b=0$ if $k=1$ -- because $M_{b}(t_0)=M_{b}(t_k)=0$ (this is the
only place where the choice of $\{ m_k \}$ really makes the difference),
and
\bee
\ba{lll}
\ds  f_{b}^{(2)}(t^*;t^{k-1};l^{k-1}) &=& \ds -\int_{T_0}^{t_1}
 F_1(t^*;s,t^{k-1};q_k,l^{k-1})M_b(s)ds \\
&-&\ds \sum_{j=2}^{k-1} \int_{t_{j-1}}^{t_{j}}
 F_j(t^*;\ldots,t_{j-1},s,t_{j},\ldots;l^{k-1,j})M_b(s)ds\\
&-&\ds \int_{t_{k-1}}^{t_{k}} F_k(t^*;t^{k-1},s;l^{k-1},q_k)M_b(s)ds.
\ea
\bee
Note that if the operators $B_l$ commute with each other, then
$f_b^{(1)}(t^*;t^{k-1};l^{k-1})$ is identically equal to zero  for all $k$.

Since $ |\a(i_{|\a|},q_{|\a|})|=|\a|-1$ and
 $\a! \geq \a(i_{|\a|},q_{|\a|})!$, it now follows from (\ref{pf:25}) that
\bee
\ba{l}
\ds \sum_{|\a|=k;i_k^{\a}=b}\frac{|\f_{\a}(t^*)|^2}{\a!}  \\ =
\ds \sum_{|\a|=k;i_k^{\a}=b}\sum_{q_k=1}^r \Big{|} \frac{1}{\sqrt{\a !}}\
\sum_{l^{k-1}}\int_{T_0}^{(k-1,t^*)}
(f^{(1)}_b+f^{(2)}_b)E_{\a(b,q_k)}dt^{k-1} \Big{|}^2  \\ \leq
\ds \sum_{q_k=1}^r \sum_{|\beta|=k-1} \Big{|} \frac{1}{\sqrt{\beta !}}\
\sum_{l^{k-1}}\int_{T_0}^{(k-1,t^*)}
 (f^{(1)}_b+f^{(2)}_b) E_{\beta}dt^{k-1} \Big{|}^2,
\ea
\bee
\ni and the proof of Proposition \ref{S.th} shows that
the last expression is equal to
\be
\label{pf:f}
\sum_{q_k=1}^r\sum_{l^{k-1}}\int_{T_0}^{(k-1,t^*)}
\Big{|} f^{(1)}_b(t^*; t^{k-1}; l^{k-1})+
f^{(2)}_b(t^*;t^{k-1}; l^{k-1}) \Big{|}^2 dt^{k-1}.
\ee

Definition of $f_b^{(1)}$  implies
\be
\label{pf:f11}
| f_b^{(1)}|^2 =0,\ k=1;\ \
| f_b^{(1)}|^2 \leq
\frac{k(C_2)^k\epsilon(B)(t^*-T_0)}{(b-1)^2}e^{C_1(t^*-T_0)}\;|U_0|^2,\
k \geq 2.
\ee

Next, direct computations yield
\bee
\ba{lll}
F_j(t^*;s^k;l^k) &= &
\ds \F_{t^*-s_k}B_{l_k}\ldots \F_{s_{j+1}-s_j}B_{l_j}
 A\F_{s_j-s_{j-1}}\ldots \F_{s_1-T_0}U_0\\
&-& \ds \F_{t^*-s_k}B_{l_k}\ldots A\F_{s_{j+1}-s_j} B_{l_j}
 \F_{s_j-s_{j-1}}\ldots \F_{s_1-T_0}U_0,
\ea
\bee
so that by assumption (3) of the theorem,
$|F_j(t^*;s^k;l^k)|^2 \leq C_0(C_2)^ke^{C_1(t^*-T_0)}\;|U_0|^2.$

After that the definition of $f_b^{(2)}$ implies:
\bee
\ba{lll}
\ds |f_{b}^{(2)}|^2  &\leq & \ds 4C_0k(C_2)^k e^{C_1(t^*-T_0)}|U_0|^2
(t^*-T_0) \int_{T_0}^{t^*} (M_{b}(s))^2 ds \\ &\leq&
\ds \frac{C_0k(C_2)^k (t^*-T_0)^3 }{(b-1)^2}e^{C_1(t^*-T_0)}\;|U_0|^2\;;
\ea
\bee
so, since $\ds \int\limits_{T_0}^{(k-1,t^*)} dt^{k-1}=
(t^*-T_0)^{k-1}/(k-1)!$, (\ref{pf:f}),
 (\ref{pf:f11}) and  the last inequality yield
\bee
\ba{l}
\E |U^N(t^*) - U_N^n(t^*)|^2 =
\ds \sum_{b \geq n+1}\sum_{k = 1}^N\sum_{|\a|=k;i_k^{\a}=b}
 \frac{\E|\f_{\a}(t^*)|^2}{\a!}  \\
\leq \ds C_2 r e^{C_1(t^*-T_0)}\Big[\epsilon(B)
(t^*-T_0^2)\sum_{k\geq 0}\frac{k+2}{k+1}
 \frac{(C_2r(t^*-T_0))^k}{k!}\\
+\ds C_0(t^*-T_0)^3 \sum_{k\geq 0}
 \frac{(k+1)(C_2r(t^*-T_0))^{k}}{k!}\Big]\E|U_0|^2
\sum_{b \geq n} \frac{1}{b^2} \\
\leq\ds \frac{2C_2r\;e^{\Cb(t^*-T_0)}}{n}\left[\epsilon(B)(t^*-T_0)^2
 +(1+(t^*-T_0)C_2r)C_0(t^*-T_0)^3 \right]\E|U_0|^2.
\ea
\bee
This completes the proof of Proposition \ref{prop:nN}. The statement of Theorem
\ref{utrn.cor} now follows from Propositions \ref{utrN.th} and \ref{prop:nN}.

\section*{Acknowledgement}
The author is grateful to Professor Boris Rozovskii for very helpful
discussions. The work was partially
 supported by the NSF grant DMS-9972016



\end{document}